\documentclass[15pt]{article}
\usepackage{amsmath, amsthm, amssymb}
\usepackage{amsmath}
\usepackage{amssymb}

\def\benr{\begin{eqnarray}}
\def\eenr{\end{eqnarray}}
\def\nn{\nonumber}

\newcommand{\norm}[2]{\|\, #1 \,\|_{#2}}

\newtheorem{proposition}{Proposition}
\newtheorem{lemma}{Lemma}
\newtheorem{theorem}{Theorem}

\begin{document}
\title{Linear bound for the dyadic paraproduct on
weighted Lebesgue space $L_2(w)$}
\author{Oleksandra V. Beznosova\footnote{Special thank to my graduate adviser Mar\'{i}a Cristina Pereyra}}

\maketitle

\begin{abstract}
The dyadic paraproduct is bounded in weighted Lebesgue spaces
$L_p(w)$ if and only if the weight $w$ belongs to the Muckenhoupt
class $A_p^d$. However, the sharp bounds on the norm of the dyadic
paraproduct are not known even in the simplest $L_2(w)$ case. In
this paper we prove the linear bound on the norm of the dyadic
paraproduct in the weighted Lebesgue space $L_2(w)$ using Bellman
function techniques and extrapolate this result to the $L_p(w)$
case.
\end{abstract}

\section{Introduction}
Let $D$ be the collection of dyadic intervals $D = \left \{I = \left
[k2^{-j}; (k+1)2^{-j}\right )\; |\; k,j \in \mathbf{Z} \right\}$,
and let $m_I f$ stand for the average of a locally integrable
function $f$ over interval $I\;$ $m_I f := \frac{1}{|I|} \int_I
{f}$.

The dyadic paraproduct is defined as
$$
\pi_b f := \sum_{I\in D}{m_If \; b_I\; h_I}
$$
where $\{ h_I\}_{I\in D}$ is the Haar basis normalized in $L_2$:
$$
h_I(x) = \frac{1}{\sqrt{|I|}} \left ( \chi_{I^+}(x) - \chi_{I^-}(x)
\right )
$$
$I^+$ and $I^-$ are left and right halves of the dyadic interval
$I$, $b_I := \langle b,h_I \rangle$ where $\langle , \rangle$ stands
for the dot product in the unweighted $L_2$, and $b$ is a locally
integrable function.

In order for the paraproduct to be bounded on $L_p$ we need $b$ to
be in $BMO^d$ i.e.:
$$
\norm{b}{BMO^d} \; := \; \left( \sup_{I} \frac{1}{|I|}
\int_{I}{\left| b(x) - m_Ib \right|^2 dx} \right)^{1/2} \; < \;
\infty.
$$
We are going to use the fact that the $BMO^d$ norm of $b$ can also
be written as:
$$
\norm{b}{BMO^d}^2 \; = \; \sup_{J\in D} \frac{1}{|J|} \sum_{I\in
D(J)}{b_I^2}.
$$

Paraproducts first appeared in the work of Bony in relation with
nonlinear partial differential equations (see \cite{Bony:81}) and
since then took one of the central places in harmonic analysis. Due
to the celebrated T(1) theorem of David and Journ\'{e}
\cite{DavidJourne:84} a singular integral operator $T$ can be
written as $T = L + \pi_{b_1} + \pi^*_{b_2}$ where $L$ is almost
translation invariant (convolution) operator, ($L1=0=L^*1$), $b_1$
is the value of $T$ at $1$ and $b_2 = T^*(1)$. The dyadic version of
this theorem can be found in \cite{Pereyra:01}. So, if one is
looking for a bound on the norms of some reasonably large class of
singular integral operators it is natural to start with the
paraproduct and with its simple dyadic "toy" model.

In this paper we prove the linear bound on the norm of dyadic
paraproduct on the weighted spaces $L_2(w)$ in terms of the $A_2^d$
characterization of the weight $w$. And now in order to prove the
linear bounds on the norms of operators with standard kernels in the
dyadic case one has to concentrate on the operator $L$.

Paraproduct holds the key to the class of singular integrals with
standard kernels. A typical representative of which is the Hilbert
transform defined by
$$
Hf(x) = P.V.\frac{1}{\pi} \int{\frac{f(x)}{x-y} dy}.
$$

Helson \& Szeg\"{o} in \cite{HelsonSzego:60} gave necessary and
sufficient condition for a weight $w$ so that $H$ maps $L^2(w)$ into
itself continuously.

In 1973, Hunt, Muckenhoupt and Wheeden (see
\cite{HuntMuckenhouptWheeden:73}) presented a new proof, where for
the first time the $A_p$-condition for the weights appeared as
necessary and sufficient condition for the boundedness of the
Hilbert transform in $L_p(w)$
$$
w\in A_p^d \; \Leftrightarrow \; \norm{w}{A_p} := \sup_{I\in D}
\left ( \frac{1}{|I|}\int_{I}{w} \right ) \left ( \frac{1}{|I|}
\int_{I}{w^{- \frac{1}{p-1} }} \right )^{p-1} \; < \infty .
$$

And a year after in \cite{CoifmanFefferman:74} Coifman and Fefferman
extended this result to a larger class of operators.

The question that has been asked is:

How is the norm of a singular operator in the weighted Lebesque
spaces $L_p(w)$ related to the Muckenhoupt $(A_p)$ characteristic of
the weight $w$, $\norm{w}{A_p^d}$. More precisely, what we need is
the sharp function $\varphi(x)$ in terms of the growth, such that
$$
\norm{Tf}{L_p(w)} \; \leq \; C \varphi \left( \norm{w}{A_p^d}
\right) \norm{f}{L_p(w)}.
$$

This kind of estimates for different singular operators is used a
lot in partial differential equations, see
\cite{FeffermanKenigPipher:91}, \cite{AstalaIvanecSaksman:01},
\cite{PetermichlVolberg:02}, \cite{BanuelosJanakiraman} and
\cite{DragicevicPetermichlVolberg:06}. Some partial answers were
given to this question.

For the Hilbert transform, Buckley showed power 2 in
\cite{Buckley:93}, Petermichl and Pott in \cite{PetermichlPott:02}
improved the exponent of $\norm{w}{A_2}$ from $2$ to $\frac{3}{2}$
and in 2006 Petermichl got the sharp power $1$ for the Hilbert
transform, see \cite{PetermichlHilbert:06}.

Later in \cite{PetermichlRiesz:06} Petermichl used similar ideas to
show linear bound for the norm of the Riesz transforms.

It was also shown that the norm of the Martingale transform on the
weighted space $L_2(w)$ depends linearly on the $\norm{w}{A_2}$, see
\cite{Wittwer:00}.

So now we can claim that singular integral operators related to the
above transforms via T(1) theorem admit linear bounds on their
norms, i.e. if $T-\pi_{b_1}-\pi_{b_2}^*$ is good enough (one of the
operators, for which we know the bound is linear), then

$$
\norm{T}{L_2(w) \rightarrow L_2(w)} \leq C \norm{w}{A_2^d}
$$

Boundedness of the paraproduct on the weighted $L_p(w)$ has been
known for a long time, a direct proof of it can be found, for
example, in \cite{KatzPereyra:97}. The best known bound on the norm
of the dyadic paraproduct so far is

$$
\norm{\pi_b}{L_2(w)\rightarrow L_2(w)} \leq C \phi(\norm{w}{A_2})
\norm{b}{BMO^d}
$$
with $\phi (x) = x^2$ and it can be found in
\cite{DragicevicGrafakosPereyraPetermichl:05}.

First we were able to improve the above result from $\phi(x) = x^2$
to $\phi(x) = x^{3/2}$ without making any significant changes to the
structure of the proof. Then using the suggestion of F.~Nazarov we
tried the duality approach which allowed us to recover $3/2$ in
multiple ways and using the version of the bilinear embedding
theorem from \cite{PetermichlHilbert:06} we were able to improve to
$\phi(x)=x(1+\log^{1/2}x)$. Using the sharp version of the bilinear
embedding theorem from \cite{NazarovTreilVolberg:99} slightly
improved the power of the logarithm in the bound
($\phi(x)=x(1+\log^{1/4}x)$). And finally, the theorem presented in
this paper shows the linear bound and in fact can rely on either one
of the bilinear embedding theorems, the one by Nazarov, Treil and
Volberg or the one from Petermichl's paper. We would also like to thank 
S. Treil for a useful conversation. 

Let us state the main result now.

\begin{theorem} (Main result)
\label{Main result theorem} The norm of dyadic paraproduct on the
weighted Lebesgue space $L_2(w)$ is bounded from above by a constant
multiple of the $A_2^d$ characteristic of the weight $w$ times the
$BMO^d$ norm of $b$, i.e. for all $f\in L_2(w)$ and all $g\in
L_2(w^{-1})$
\begin{eqnarray}
\label{Para dual linear} \left \langle \pi_b f, g \right
\rangle_{L_2} \leq C \;\norm{w}{A_2} \;\norm{b}{BMO^d}\;
\norm{f}{L_2(w)} \;\norm{g}{L_2(w^{-1})}.
\end{eqnarray}
\end{theorem}

Which together with the sharp version of the Rubio De Francia's
extrapolation theorem from
\cite{DragicevicGrafakosPereyraPetermichl:05} produces $L_p$ bounds
of the following type:

\begin{theorem}
\label{Lp bound theorem} Let $w\in A_p^d$ and $b\in BMO^d$. Then the
norm of dyadic paraproduct $\pi_b$ on the weighted $L_p(w)$ is
bounded by
$$
\norm{\pi_b}{L_p(w)\rightarrow L_p(w)} \; \leq \; C_1(p)
\norm{w}{A_p^d}  \norm{b}{BMO^d} \; \; \; \; when \; p\geq 2
$$
and by
$$
\norm{\pi_b}{L_p(w)\rightarrow L_p(w)} \; \leq \; C_2(p)
\norm{w}{A_p^d}^{\frac{1}{p-1}} \norm{b}{BMO^d} \; \; \; \; when \;
p<2,
$$
where $C_1(p)$ and $C_2(p)$ are constants that only depend on $p$.
\end{theorem}

This paper is constructed as follows:

Section 2: proof of the main result based on three propositions.

Section 3: Bellman function proof of Proposition \ref{My little
lemma}.

Section 4: Bellman function proof of Proposition \ref{prop.sum1}.

Section 5: Bellman function proof of Proposition \ref{prop.sum2}.

\section{Proof of the main result}

\begin{proof}
In order to prove Theorem \ref{Main result theorem} it is enough to
show that $\forall \; f,g \in L_2$
$$
\left \langle \pi_b \left ( fw^{-1/2}\right) ; gw^{1/2} \right
\rangle \; \leq \; C \norm{w}{A_2^d} \, \norm{b}{BMO^d} \,
\norm{f}{2} \, \norm{g}{2},
$$
where $\left \langle \pi_b \left ( fw^{-1/2}\right) ; gw^{1/2}
\right \rangle$ can be written as the following sum
$$
\left \langle \pi_b \left ( fw^{-1/2}\right) ; gw^{1/2} \right
\rangle \; = \; \sum_{I\in D}{m_I \left( fw^{-1/2} \right) b_I \left
\langle gw^{1/2};h_I \right \rangle} \; =: \; {\sum}_1
$$

Now, we are going to decompose this sum into parts using weighted
Haar system of functions:

Let $H_I^w$ be defined in the following way:
$$
H_I^w \; := \; h_I \sqrt{|I|} - A_I^w \chi_I.
$$
In order to make $\{H_I^w\}$ an orthogonal system of functions in
$L_2(w)$, we take $A_I^w$ to be
$$
A_I^w \; := \; \frac{m_{I^+}w - m_{I^-}w}{2m_Iw},
$$
then $\{w^{1/2} H_I^w\}$ is orthogonal in $L_2$ with norms bounded
from above by $\norm{w^{1/2} H_I^w}{L_2} \leq \sqrt{|I| m_Iw}$.

Then by Bessel's inequality we have:

\begin{equation}\label{Bessel's}
\forall g\in L_2 \; \;\;\; \; \sum_{I\in D}{\frac{1}{|I|m_Iw}
\left\langle g; w^{1/2}H_I^w \right\rangle_{L_2}^2} \; \leq \;
\norm{g}{L_2}^2.
\end{equation}

The weighted Haar functions were first introduced in
\cite{CoifmanJonesSemmes:89} and are extremely useful in weighted
inequalities, see \cite{NazarovTreilVolberg:99} and
\cite{Pereyra:05}.

We can break ${\sum}_1$ into two sums:
\begin{eqnarray*}
{\sum}_1 \; &=& \; \sum_{I\in D}{m_I \left( fw^{-1/2} \right) b_I
\left \langle gw^{1/2};h_I \right \rangle} \\
&=& \sum_{I\in D}{m_I \left( fw^{-1/2} \right) b_I
\frac{1}{\sqrt{|I|}} \left \langle g;w^{1/2}H_I^w \right \rangle} +
\sum_{I\in D}{m_I \left( fw^{-1/2} \right) b_I \frac{1}{\sqrt{|I|}}
\left \langle gw^{1/2};A_I^w \chi_I \right \rangle}  \\
&=:& \; {\sum}_2 + {\sum}_3.
\end{eqnarray*}

And now we will bound ${\sum}_2$ and ${\sum}_3$.

We claim that both sums, ${\sum}_2$ and ${\sum}_3$, depend on the
$\norm{w}{A_2^d}$ at most linearly:
\begin{equation}
\label{main proof sum1} {\sum}_2 \; = \; \sum_{I\in D}{m_I \left(
fw^{-1/2} \right) b_I \frac{1}{\sqrt{|I|}} \left \langle
g;w^{1/2}H_I^w \right \rangle} \; \leq \; C\norm{w}{A_2^d}
\norm{b}{BMO^d} \norm{f}{L_2} \norm{g}{L_2}.
\end{equation}
and
\begin{equation}
\label{main proof sum2} {\sum}_3 = \sum_{I\in D}{m_I \left(
fw^{-1/2} \right) b_I A_I^w \sqrt{|I|} m_I \left( gw^{1/2} \right)}
\leq C \norm{w}{A_2^d} \norm{b}{BMO^d} \norm{f}{L_2} \norm{g}{L_2}.
\end{equation}
Before going into the proofs of (\ref{main proof sum1}) and
(\ref{main proof sum2}) let us analyze the partition ${\sum}_1 =
{\sum}_2 + {\sum}_3$.

The sum ${\sum}_2$ is close to the "weighted" version of a
paraproduct over a weighted space $L_2(w)$, which behaves similar to
the unweighted situation, while ${\sum}_3$ takes into account the
difference between the norm of the paraproduct on weighted and
unweighted $L_2$. In the simplest case $w=const$, $\norm{w}{A_2^d} =
1$, ${\sum}_1 = {\sum}_2$, and we recover classical results, while
${\sum}_3 = 0$.

Note also, that for weights with small $A_2^d$-characteristics
${\sum}_2$ will be dominating and ${\sum}_3$ will be close to $0$,
while for $\norm{w}{A_2^d}$ large ${\sum}_3$ becomes more important.

Bound on ${\sum}_2$ is very straight-forward and very similar to the
classical case. We decompose ${\sum}_2$ into the product of two sums
using Cauchy-Schwarz:
\begin{eqnarray*}
{\sum}_2 \; &=& \; \sum_{I\in D}{m_I \left( fw^{-1/2} \right) b_I
\frac{1}{\sqrt{|I|}} \left \langle g;w^{1/2}H_I^w \right \rangle}
\\
&\leq& \; \left( \sum_{I\in D}{m_I^2\left( fw^{-1/2} \right)
b_I^2 m_Iw} \right )^{1/2} \left( \sum_{I\in D}{\frac{1}{|I|m_Iw}
\left\langle g; w^{1/2}H_I^w \right\rangle^2} \right)^{1/2}.
\end{eqnarray*}
By (\ref{Bessel's})
$$
\sum_{I\in D}{\frac{1}{|I|m_Iw} \left\langle g; w^{1/2}H_I^w
\right\rangle^2} \; \leq \; \norm{g}{L_2}^2.
$$
So, for (\ref{main proof sum1}) it is enough to show that
\begin{equation}
\label{main proof sum3} \sum_{I\in D}{m_I^2\left( fw^{-1/2} \right)
b_I^2 m_Iw} \; \leq \; C \norm{w}{A_2^d}^2 \norm{b}{BMO^d}^2
\norm{f}{L_2}^2.
\end{equation}
By the weighted Carleson embedding theorem, which can be found, for
example, in \cite{NazarovTreilVolberg:99}, and $(d\sigma)$ version
of it can be found in \cite{Pereyra:05}, (\ref{main proof sum3})
holds if and only if
$$
\forall J\in D \; \; \; \; \frac{1}{|J|} \sum_{I\in D(J)}
{m_I^2w^{-1} \; m_Iw \; b_I^2} \; \leq \; C \norm{w}{A_2^d}^2
\norm{b}{BMO^d}\; m_Jw^{-1}.
$$
And since $\forall I\in D \;$ $\; m_Iw\; m_Iw^{-1} \leq
\norm{w}{A_2^d}$, it is enough to verify that
\begin{equation}
\label{main proof sum4} \forall J\in D \;\;\;\; \frac{1}{|J|}
\sum_{I\in D(J)}{m_Iw^{-1}\; b_I^2} \; \leq \; C \norm{w}{A_2^d}
\norm{b}{BMO^d}\; m_Jw^{-1}.
\end{equation}
Inequality (\ref{main proof sum4}) follows from the fact that $b\in
BMO^d$ and hence the sequence $\{ b_I^2\}_{I\in D}$ is a Carleson
sequence with Carleson constant $\norm{b}{BMO^d}^2$:
\begin{equation}
\label{Carleson for b} \forall J\in D \;\;\;\;
\frac{1}{|J|}\sum_{I\in D(J)}{b_I^2} \; \leq \; \norm{b}{BMO^d}^2,
\end{equation}
and the following proposition, which we are going to prove using the
Bellman function technique in Section 3.
\begin{proposition}
\label{My little lemma} Let $w\in A_2^d$ and $\{ \lambda_I\}$ be a
Carleson sequence of nonnegative numbers
$$
\forall J\in D \;\;\;\; \frac{1}{|J|}\sum_{I\in D(J)}{\lambda_I} \;
\leq \; Q,
$$
then $\forall J\in D$
\begin{equation}
\label{little lemma main inequality} \frac{1}{|J|} \sum_{I\in
D(J)}{\frac{\lambda_I}{m_Iw^{-1}}} \; \leq \; 4Q\; m_Jw
\end{equation}
and
\begin{equation}
\label{little lemma consequence} \frac{1}{|J|} \sum_{I\in D(J)}{m_Iw
\lambda_I} \; \leq \; 4Q\norm{w}{A_2^d}\; m_Jw.
\end{equation}
\end{proposition}

Estimate (\ref{little lemma consequence}) applied to $\lambda_I =
b_I^2$ and $w^{-1}$ ($w^{-1} \in A_2^d$ and $\norm{w^{-1}}{A_2^d} =
\norm{w}{A_2^d}$) provides (\ref{main proof sum4}), so bound
(\ref{main proof sum1}) on ${\sum}_2$ holds.

Now we need to prove bound (\ref{main proof sum2}) on ${\sum}_3$. It
is a little bit more involved. We want to show that
$$
{\sum}_3 = \sum_{I\in D}{b_I A_I^w \sqrt{|I|}\; m_I\left( fw^{-1/2}
\right) \; m_I\left( gw^{1/2}\right)} \leq C\norm{w}{A_2^d}
\norm{b}{BMO^d} \norm{f}{2} \norm{g}{2}.
$$

We are going to use a version of the bilinear embedding theorem by
Petermichl from \cite{PetermichlHilbert:06}:

\begin{theorem} (Petermichl) Let $w\in A_2$, $\norm{w}{A_2} \leq Q$. Let
$\{ \alpha_I\}_{I\in D}$ be a sequence of non-negative numbers such
that:
$$
\forall J\in D \;\;\; \frac{1}{|J|}\sum_{I\in D(J)}{\alpha_I \; m_Iw
\; m_Iw^{-1}} \; \leq \; Q,
$$
$$
\forall J\in D \;\;\; \frac{1}{|J|}\sum_{I\in D(J)}{\alpha_I \; m_Iw
} \; \leq \; Q m_Jw,
$$
$$
\forall J\in D \;\;\; \frac{1}{|J|}\sum_{I\in D(J)}{\alpha_I \;
m_Iw^{-1}} \; \leq \; Q m_Jw^{-1},
$$
then there is a constant $C>0$ such that $\forall f,g\in L_2$
$$
\sum_{I\in D}{\alpha_I \; m_I\left( fw^{-1/2}\right) \; m_I\left(
gw^{1/2} \right)} \; \leq \; CQ\norm{f}{L_2}\norm{g}{L_2}.
$$
\end{theorem}

So, in order to complete the proof it is enough to show that the
following three bounds hold:

\begin{equation}
\label{main proof sum5} \forall J\in D \;\;\;
\frac{1}{|J|}\sum_{I\in D(J)}{|b_I A_I^w| \sqrt{|I|}\; m_Iw \;
m_Iw^{-1}} \; \leq \; C \norm{w}{A_2^d},
\end{equation}
\begin{equation}
\label{main proof sum6} \forall J\in D \;\;\;
\frac{1}{|J|}\sum_{I\in D(J)}{|b_I A_I^w| \sqrt{|I|} \; m_Iw } \;
\leq \; C \norm{w}{A_2^d} m_Jw,
\end{equation}
\begin{equation}
\label{main proof sum7} \forall J\in D \;\;\;
\frac{1}{|J|}\sum_{I\in D(J)}{|b_I A_I^w| \sqrt{|I|} \; m_Iw^{-1}}
\; \leq \; C \norm{w}{A_2^d} m_Jw^{-1},
\end{equation}

The following Proposition helps us handle the first sum (\ref{main
proof sum5}).

\begin{proposition} \label{prop.sum1}
Let $w$ be a weight from $A_2^d$, then $\forall J\in D$
\begin{equation}
\nn \frac{1}{|J|}\sum_{I\in D(J)}{\left( \frac{m_{I^+}w -
m_{I^-}w}{m_Iw} \right)^2 |I| m_I^{1/4}w \; m_I^{1/4}w^{-1}} \; \leq
\; C m_J^{1/4}w \; m_J^{1/4}w^{-1}.
\end{equation}
\end{proposition}

Note that a simple consequence of Proposition \ref{prop.sum1} is
$$
\frac{1}{|J|}\sum_{I\in D(J)}{\left( \frac{m_{I^+}w -
m_{I^-}w}{m_Iw} \right)^2 |I| m_Iw \; m_Iw^{-1}} \; \leq \; C
\norm{w}{A_2^d}^{3/4} m_J^{1/4}w \; m_J^{1/4}w^{-1}
$$
and hence
\begin{equation}
\label{consequence hence} \frac{1}{|J|}\sum_{I\in D(J)}{\left(
\frac{m_{I^+}w - m_{I^-}w}{m_Iw} \right)^2 |I| m_Iw \; m_Iw^{-1}} \;
\leq \; C \norm{w}{A_2^d}.
\end{equation}

Then by Cauchy-Schwarz
$$
\frac{1}{|J|} \sum_{I\in D(J)}{|b_I A_I^w| \sqrt{|I|} m_Iw \,
m_Iw^{-1}} \leq \left( \frac{1}{|J|} \sum_{I\in D(J)}{b_I^2 m_Iw \,
m_Iw^{-1}} \right)^{\frac{1}{2}} \left( \frac{1}{|J|} \sum_{I\in
D(J)}{(A_I^w)^2 |I| m_Iw \, m_Iw^{-1}} \right)^{\frac{1}{2}},
$$
by (\ref{consequence hence})
$$
\frac{1}{|J|} \sum_{I\in D(J)}{(A_I^w)^2 |I|\, m_Iw \; m_Iw^{-1}} \;
\leq \; C \norm{w}{A_2^d},
$$
and by (\ref{Carleson for b})
$$
\frac{1}{|J|} \sum_{I\in D(J)}{b_I^2 \, m_Iw \; m_Iw^{-1}} \; \leq
\; \norm{w}{A_2^d} \frac{1}{|J|} \sum_{I\in D(J)}{b_I^2} \; \leq \;
\norm{w}{A_2^d} \, \norm{b}{BMO^d}^2.
$$

Linear bound on the second sum (\ref{main proof sum6}) follows by
Cauchy-Schwarz, from the sharp result by J.Wittwer
\cite{Wittwer:00}:

\begin{lemma} (J.Wittwer)
Let $w\in A^d_2$ be a weight, then
$$
\forall J\in D \;\;\; \frac{1}{|J|}\sum_{I\in
D(J)}{\left(\frac{m_{I_-}w - m_{I_+}w}{2m_Iw} \right)^2 |I|}\; m_I w
\; \leq \;  C \norm{w}{A_2^d} m_J w
$$
and this result is sharp.
\end{lemma}
\noindent and Proposition \ref{My little lemma} (inequality
(\ref{little lemma consequence})) applied to $\lambda_I = b_I^2$:
$$
\forall J\in D \;\;\;\; \frac{1}{|J|} \sum_{I\in D(J)}{b_I^2 \;
m_Iw} \; \leq \; C \norm{w}{A_2^d} \; \norm{b}{BMO^d}^2 \; m_Jw.
$$

And the next proposition together with (\ref{main proof sum4})
allows us to bound the third sum (\ref{main proof sum7}) in a
similar way.

\begin{proposition}
\label{prop.sum2} Let $w$ be a weight in $A_2^d$, then for all
dyadic intervals $J$:
$$
\frac{1}{|J|}\sum_{I\in D(J)}{\left( \frac{m_{I^+}w -
m_{I^-}w}{m_Iw} \right)^2 |I| m_Iw^{-1}} \; \leq \; C
\norm{w}{A_2^d} m_Jw^{-1}.
$$
\end{proposition}

Which completes the proof of the Theorem \ref{Main result theorem}.
\end{proof}

\section{Bellman function proof of Proposition \ref{My little lemma}}

We are going to show that for any Carleson sequence $\{\lambda_I
\}_{I\in D}$ with constant $Q$, $\lambda_I \geq 0$
$$
\forall J\in D \;\;\;\; \frac{1}{|J|} \sum_{I\in D(J)}{\lambda_I} \;
\leq \; Q
$$
the inequality (\ref{little lemma main inequality}) holds for any
dyadic interval $J$:
$$
\frac{1}{|J|} \sum_{I\in D(J)}{\frac{\lambda_I}{m_Iw^{-1}}} \; \leq
\; 4Qm_Jw.
$$
Note that inequality (\ref{little lemma consequence}) follows from
inequality (\ref{little lemma main inequality}).
\begin{lemma}
\label{prop1 lemma1} Suppose there exists a real valued function of
3 variables $B(x) = B(u,v,l)$, whose domain $\mathcal{D}$ is given
by those $x=(u,v,l)\in \mathbb{R}^3$ such that
$$
u,\; v,\;l\; \geq \; 0,
$$
$$
uv \; \geq \; 1,
$$
$$
l \; \leq \; 1,
$$
whose range is given by
$$
0\; \leq \; B(x) \; \leq \; m_Jw,
$$
and such that the following convexity property holds:
$$
\forall x,\;x_\pm \in \mathcal{D} \;\;\; such \;\; that \;\;\; x -
\frac{x_+ + x_-}{2} = (0,0,\alpha)
$$
\begin{equation}
\label{prop1 lemma1 convex cond} B(x) - \frac{B(x_+) + B(x_-)}{2} \;
\geq \; \frac{1}{4v} \alpha
\end{equation}
Then Proposition \ref{My little lemma} holds.
\end{lemma}

\begin{proof}[Proof of Lemma \ref{prop1 lemma1}]

Fix a dyadic interval $J$. Let $x_J = (u_J, v_J, l_J)$ where $u_J =
m_Jw$, $v_J = m_Jw^{-1}$ and $l_J=\frac{1}{|J|Q}\sum_{I\in
D(J)}{\lambda_I}$. Clearly for each dyadic $J$, $x_J$ belongs to the
domain $\mathcal{D}$. Let $x_\pm := x_{J^\pm} \in \mathcal{D}$. By
definition,
$$
x_J - \frac{x_{J^+} + x_{J^-}}{2} \; = \; (0,0,\alpha_J),
$$
where $\alpha_J := \frac{1}{|J|Q}\lambda_J$. Then, by convexity
condition (\ref{prop1 lemma1 convex cond})
\begin{eqnarray*}
m_Jw \; \geq \; B(x_J) &\geq& \frac{B(x_{J^+})}{2} + \frac{B(x_{J^-})}{2}
+ \frac{1}{4v_J}\alpha_J\\
&=& \frac{B(x_{J^+})}{2} + \frac{B(x_{J^-})}{2} +
\frac{1}{4|J|Qm_Jw^{-1}}\lambda_J.
\end{eqnarray*}
Iterating this procedure and using the assumption that $B\geq 0$ on
$\mathcal{D}$ we get:
$$
m_Jw \geq \frac{1}{4|J|Q}\sum_{I\in
D(J)}{\frac{\lambda_I}{m_Iw^{-1}}}
$$
which implies Proposition \ref{My little lemma}.
\end{proof}
So, Proposition \ref{My little lemma} will hold if we can show
existence of the function $B$ of the Bellman type, satisfying the
conditions of Lemma \ref{prop1 lemma1}.

\begin{lemma}
\label{prop1 lemma2} The following function
$$
B(u,v,l) := u - \frac{1}{v(1+l)}
$$
is defined on $\mathcal{D}$, $0 \leq B(x) \leq u$ for all $x=(u,v,l)
\in \mathcal{D}$ and satisfies the following differential
inequalities on $\mathcal{D}$:
\begin{equation}
\label{prop1 lemma2 ineq1} \frac{\partial B}{\partial l} \; \geq \;
\frac{1}{4v}
\end{equation}
and
\begin{equation}
\label{prop1 lemma2 ineq2} -d^2B \; \geq \; 0.
\end{equation}
Moreover, conditions (\ref{prop1 lemma2 ineq1}) and (\ref{prop1
lemma2 ineq2}) imply the convexity condition (\ref{prop1 lemma1
convex cond}).
\end{lemma}

\begin{proof}
Range conditions are easy to verify: since all variables are
positive on $\mathcal{D}$ and $uv\geq 1$, we have
$$
0\; \leq \; B(u,v,l)=\frac{uv(1+l) - 1}{v(1+l)} = u -
\frac{1}{v(1+l)} \; \leq \; u.
$$
It is nothing but a calculus exercise to check the differential
conditions:
$$
\frac{\partial B}{\partial l} \; = \; \frac{1}{v(1+l)^2} \; \geq \;
\frac{1}{4v}
$$
since $l\geq 1$. And
\begin{equation}
\nn
 -d^2B \; = \; (du, dv, dl) \left(
\begin{tabular}{ c c c }
$0$ & $0$ & $0$ \\
$0$ & $\frac{2}{v^3(1+l)}$ & $\frac{1}{v^2(1+l)^2}$ \\
$0$ & $\frac{1}{v^2(1+l)^2}$ & $\frac{2}{v(1+l)^3}$ \\
\end{tabular}
 \right) \left(
 \begin{tabular}{ c }
 $du$ \\
 $dv$ \\
 $dl$ \\
 \end{tabular} \right) \; \geq \; 0
\end{equation}

And finally let us see how differential conditions (\ref{prop1
lemma2 ineq1}) and (\ref{prop1 lemma2 ineq2}) imply the convexity
condition (\ref{prop1 lemma1 convex cond}):
$$
B(x) - \frac{B(x_+) + B(x_-)}{2} = \left[ B(x) - B(\frac{x_+ +
x_-}{2})\right] + \left[ B(\frac{x_+ + x_-}{2}) - \frac{B(x_+) +
B(x_-)}{2} \right] =
$$
$$
= \frac{\partial B}{\partial l}(u,v,l^\prime)\alpha -
\int_{-1}^1{(1-|t|)b^{\prime\prime}(t)dt},
$$
where $b(t) := B(s(t))$, $s(t) := \frac{1+t}{2}s_+ +
\frac{1-t}{2}s_-$, $-1\leq t \leq 1$, note that $s(t) \in
\mathcal{D}$ whenever $s_+$ and $s_-$ do since $\mathcal{D}$ is a
convex domain. Then differential inequalities trivially imply that
$-b^{\prime\prime}(t) \geq 0$ and
$$
B(x) - \frac{B(x_+) + B(x_-)}{2} = \frac{\partial B}{\partial
l}(u,v,l^\prime)\alpha - \int_{-1}^1{(1-|t|)b^{\prime\prime}(t)dt}
\geq \frac{1}{4v}\alpha.
$$

And proofs of both Lemma \ref{prop1 lemma2} and Proposition \ref{My
little lemma} are complete.
\end{proof}

\section{Proof of the Proposition \ref{prop.sum1}}

We are going to prove that there is a numerical constant $C>0$, such
that for all dyadic intervals $J\in D$

\begin{equation}
\label{essential part for sum1} \frac{1}{|J|}\sum_{I\in D(J)}{\left(
\frac{m_{I^+}w - m_{I^-}w}{m_Iw} \right)^2 |I| m_I^{1/4}w \;
m_I^{1/4}w^{-1}} \; \leq \; C m_J^{1/4}w \; m_J^{1/4}w^{-1}.
\end{equation}
using Bellman function technique.

\begin{lemma}
\label{Sum1 Lemma1} Suppose there exists a real-valued function of
two variables $B(x) = B(u,v)$, whose domain $\mathcal{D}$ is given
by those $x=(u,v) \in \mathbb{R}^2$ such that
\begin{equation}
\label{Sum1 domain cond 1} u,v\; \geq \; 0
\end{equation}
\begin{equation}
\label{Sum1 domain cond 2} uv\; \geq \; 1,
\end{equation}
whose range is given by
$$
0 \; \leq \; B(x) \; \leq \; \sqrt[4]{uv}, \;\;\;\; x\in
\mathcal{D},
$$
and such that the following convexity property holds:
\begin{equation}
\label{Sum1 convex cond} if \;\; x = \frac{x_+ + x_-}{2} \;\;\;\;
then \;\;\;\; B(x) - \frac{B(x_+) + B(x_-)}{2} \; \geq \; C
\frac{v^{1/4}}{u^{7/4}} (u_+ - u_-)^2
\end{equation}
with a numerical constant $C$ independent of everything, then the
Proposition \ref{prop.sum1} will be proved.
\end{lemma}

\begin{proof}
Let $u_I := m_Iw$, $v_I := m_Iw^{-1}$, $v_+ = v_{I^+}$, $v_- =
v_{I^-}$ and similarly for $u_\pm$. Then by H\"{o}lder's inequality
$(u,v)$ and $(u_\pm,v_\pm)$ belong to the $\mathcal{D}$.

Fix $J\in D$, by the convexity and range conditions
\begin{eqnarray*}
|J| \sqrt[4]{m_Jw \; m_Jw^{-1}} \; &\geq& \; |J| B(u_J, v_J) \\
&\geq& \; \frac{|J|}{2}B(u_+, v_+) + \frac{|J|}{2}B(u_-, v_-) + |J|
C \frac{m_J^{1/4}w^{-1}}{m_J^{7/4}w} \left( m_{J^+}w - m_{J^-}w
\right)^2 \\
&=& \; |J^+|B(u_+, v_+) + |J^-|B(u_-, v_-) + |J| C
\frac{m_J^{1/4}w^{-1}}{m_J^{7/4}w} \left( m_{J^+}w - m_{J^-}w
\right)^2.
\end{eqnarray*}
Iterating this process and using the fact that $B(u,v) \geq 0$ we
get:
$$
|J| \sqrt[4]{m_Jw \; m_Jw^{-1}} \; \geq \; C \sum_{I\in
D(J)}{|I|\frac{m_I^{1/4}w^{-1}}{m_I^{7/4}w} \left( m_{J^+}w -
m_{J^-}w \right)^2},
$$
which completes the proof of Lemma \ref{Sum1 Lemma1}.
\end{proof}

Now, in order to complete the proof of (\ref{essential part for
sum1}) we need to show existence of the Bellman type function $B$
which satisfies the conditions of Lemma \ref{Sum1 Lemma1}.

\begin{lemma}
\label{Sum1 Lemma2} The following function
$$
B(u,v) := \sqrt[4]{uv}
$$
is defined on $\mathcal{D}$, $0\leq B(u,v) \leq \sqrt[4]{uv}$ for
all $(u,v)\in \mathcal{D}$, and satisfies the following differential
inequality in $\mathcal{D}$:
\begin{equation}
\label{Sum1 diff convex cond} -d^2B \; \geq \; \frac{1}{8}
\frac{v^{1/4}}{u^{7/4}} |du|^2.
\end{equation}
Furthermore, this implies the convexity condition (\ref{Sum1 convex
cond}) of Lemma \ref{Sum1 Lemma1}.
\end{lemma}

\begin{proof}

Since $u$ and $v$ are positive in the domain $\mathcal{D}$, function
$B=\sqrt[4]{uv}$ is well defined on $\mathcal{D}$ and condition
$0\leq B(u,v) \leq \sqrt[4]{uv}$ is trivially satisfied.

Let us prove the differential inequality (\ref{Sum1 diff convex
cond}) now:
\begin{equation}
\nn
 -d^2B \; = \; \frac{1}{16} (du, dv) \left(
\begin{tabular}{ c c }
$3v^{\frac{1}{4}} u^{\frac{-7}{4}}$ & $-v^{\frac{-3}{4}} u^{\frac{-3}{4}}$ \\
$-v^{\frac{-3}{4}}u^{\frac{-3}{4}}$ & $3v^{\frac{-7}{4}}u^{\frac{1}{4}}$ \\
\end{tabular}
 \right) \left(
 \begin{tabular}{ c }
 $du$ \\
 $dv$ \\
 \end{tabular} \right)
\end{equation}
\begin{equation}
\nn = \frac{1}{8} (du, dv) \left(
\begin{tabular}{ c c }
$v^{\frac{1}{4}} u^{\frac{-7}{4}}$ & $0$ \\
$0$ & $v^{\frac{-7}{4}}u^{\frac{1}{4}}$ \\
\end{tabular}
 \right) \left(
 \begin{tabular}{ c }
 $du$ \\
 $dv$ \\
 \end{tabular} \right) + \frac{1}{16} (du, dv) \left(
\begin{tabular}{ c c }
$v^{\frac{1}{4}} u^{\frac{-7}{4}}$ & $-v^{\frac{-3}{4}} u^{\frac{-3}{4}}$ \\
$-v^{\frac{-3}{4}}u^{\frac{-3}{4}}$ & $v^{\frac{-7}{4}}u^{\frac{1}{4}}$ \\
\end{tabular}
 \right) \left(
 \begin{tabular}{ c }
 $du$ \\
 $dv$ \\
 \end{tabular} \right)
\end{equation}
$$
\geq \; \frac{1}{8} v^{\frac{1}{4}} u^{\frac{-7}{4}} |du|^2,
$$
as we wanted to show.

Now we only need to check the convexity condition (\ref{Sum1 convex
cond}). We fix an interval $I$ and let
$$
b(t) \; := \; B(u_t, v_t), \;\;\;\; -1 \; \leq \; t \; \leq \; 1,
$$
where
$$
u_t \; := \;\frac{1}{2}(t+1)u_+ + \frac{1}{2}(1-t)u_-
$$
and
$$
v_t \; := \;\frac{1}{2}(t+1)v_+ + \frac{1}{2}(1-t)v_-.
$$
What we want to show is
$$
b(0) - \frac{b(1) + b(-1)}{2} \; \geq \; C \frac{v^{1/4}}{u^{7/4}}
|du|^2.
$$
It is easy to see that
$$
b(0) - \frac{1}{2}\left( b(-1) + b(1) \right) \; = \; -\frac{1}{2}
\int_{-1}^{1}{(1+|t|) b^{\prime\prime}(t)dt}.
$$
Note that
\begin{equation}
\label{temp1} -b^{\prime\prime}(t) \; \geq \; \frac{1}{32} v_t^{1/4}
u_t^{-7/4} \left( u_1 - u_{-1}\right)^2
\end{equation}
and that $\forall t\in \left[ -1/2; 1/2\right]$
$$
u_t \; = \; u_0 + \frac{1}{2}t(u_1 - u_{-1}),
$$
since domain $\mathcal{D}$ is convex $u_t\in \mathcal{D}$, and
$$
|u_1 - u_{-1}| \; \leq \; |u_1| + |u_{-1}|, \;\;\;\; |t| \leq 1/2,
\;\; u_1, u_{-1} \geq 0,
$$
$$
-u_0 \; = \; -\frac{1}{2} (u_1 + u_{-1}) \; \leq \; t(u_1 - u_{-1})
\; \leq \; \frac{1}{2}(u_1 + u_{-1}) \; = \; u_0,
$$
so $u_t \; \leq \; \frac{3}{2}u_0$ and similarly $v_t \; \geq \;
\frac{1}{2} v_0$ for $t\in [-1/2; 1/2]$. Together with (\ref{temp1})
it makes
$$
-b^{\prime\prime}(t) \; \geq \; C v_0^{1/4} u_0^{-7/4} (u_1 -
u_{-1})^2.
$$
So,
$$
B(u,v) - \frac{1}{2}(B(u_+, v_+) - B(u_-,v_-)) \; = \; b(0) -
\frac{1}{2}(b(1)+b(-1)) \; \geq \; C\frac{v^{1/4}}{u^{7/4}}|du|^2
$$
with numerical constant $C$ independent of everything. Which
completes the proof of Lemma \ref{Sum1 Lemma2} and Proposition
\ref{prop.sum1}.
\end{proof}

\section{Proof of the Proposition \ref{prop.sum2}}

First note that since for every dyadic interval $I$ we have $m_Iw\;
m_Iw^{-1} \leq \norm{w}{A_2^d}$, it is enough to show that
\begin{equation}
\label{essential part of sum2} \forall J\in D \;\;\;\; \frac{1}{|J|}
\sum_{I\in D(J)}{\frac{(m_{I^+}w - m_{I^-}w)^2}{m_I^3w} |I|} \; \leq
\; C m_Jw^{-1}
\end{equation}
for some numerical constant $C$.
\begin{lemma}
\label{Sum2 Lemma1} Suppose there exists a real-valued function of
two variables $B(x) = B(u,v)$, whose domain $\mathcal{D}$ is given
by those $x = (u,v) \in \mathbb{R}^2$ such that
\begin{equation}
\label{Sum2 domain cond 1} u,v \; \geq \; 0,
\end{equation}
\begin{equation}
\label{Sum2 domain cond 2} uv \; \geq \; 1,
\end{equation}
whose range is given by
$$
0\; \leq \; B(x) \; \leq \; v
$$
and such that the following convexity property holds:
\begin{equation}
\label{Sum2 convex cond} if \;\;\; x = \frac{x_+ + x_-}{2} \;\;\;
then \;\;\; B(x) - \frac{B(x_+) + B(x_-)}{2}\; \geq \; C
\frac{1}{u^3}(u_+ - u_-)^2
\end{equation}
with some numerical constant $C$ independent of everything. Then
Proposition \ref{prop.sum2} will be proved (inequality
(\ref{essential part of sum2}) holds for all dyadic intervals $J$).
\end{lemma}

\begin{proof}
Let $u_I := m_Iw$, $v_I := m_Iw^{-1}$, $v_+ = v_{I^+}$, $v_- =
v_{I^-}$ and similarly for $u_\pm$. Then by H\"{o}lder's inequality
$(u,v)$ and $(u_\pm,v_\pm)$ belong to the $\mathcal{D}$.

Fix $J\in D$, by the convexity property and range conditions
\begin{eqnarray*}
|J|m_Jw^{-1} \; &\geq& \; |J|B(u_J, v_J) \\
&\geq& \; \frac{|J|}{2} B(u_+, v_+) + \frac{|J|}{2}B(u_-, v_-) +
C|J|\frac{1}{m_J^3w}(m_{J^+}w - m_{J^-}w)^2 \\
&=& \; |J^+|B(u_+, v_+) + |J^-|B(u_-, v_-) + C|J|\frac{1}{m_J^3w}
(m_{J^+}w - m_{J^-}w)^2.
\end{eqnarray*}
Iterating this process and using positivity of function $B$, we get
$$
|J|m_Jw^{-1} \; \geq \; C\sum_{I\in
D(J)}{|I|\frac{1}{m_I^3w}(m_{I^+}w - m_{I^-}w)^2},
$$
which completes the proof of Lemma \ref{Sum2 Lemma1}.
\end{proof}

To prove inequality (\ref{essential part of sum2}) and Proposition
\ref{prop.sum2} we need to show the existence of the function $B$ of
the Bellman type satisfying conditions of Lemma \ref{Sum2 Lemma1}.

\begin{lemma}
\label{Sum2 Lemma2} The following function
$$
B(u,v) \; = \; v - \frac{1}{u}
$$
defined on domain $\mathcal{D}$, $0\leq B(u,v)\leq v$ for all
$(u,v)\in \mathcal{D}$ and satisfies the following differential
inequality in $\mathcal{D}$:
$$
-d^2B \; \geq \; \frac{2}{u^3}|du|^2.
$$
Moreover, it implies the convexity condition (\ref{Sum2 convex
cond}) with some numerical constant $C$ independent of everything.
\end{lemma}

\begin{proof}
First note that since $uv\geq 1$ and $u$ and $v$ are both positive
in the domain $\mathcal{D}$, $B$ is well-defined and
$$
0\; \leq \; B(u,v) \; = \; \frac{uv-1}{u} \; = \; v-\frac{1}{u} \;
\leq \; v
$$
on $D$ and $-d^2B \; = \; 2u^{-3}|du|^2$.

Convexity condition (\ref{Sum2 convex cond}) follows from this in
practically the same way as in Proposition \ref{prop.sum1}.
\end{proof}

\bibliographystyle{plain}   

\begin{verbatim}
Oleksandra V. Beznosova
Department of Mathematics and Statistics
University of New Mexico
Albuquerque, NM 87131, USA
alexbeznosova@yahoo.com
\end{verbatim}

\end{document}